\documentclass[11pt,a4paper]{article}
\usepackage{mathrsfs}
\usepackage{mathrsfs}
\usepackage{mathrsfs}
\usepackage{amsfonts}
\usepackage{amsfonts}
\usepackage{amsfonts}
\usepackage{mathrsfs}
\usepackage{amssymb}
\usepackage{color}
\usepackage{amsmath}
\newtheorem{theorem}{Theorem}[section]

\numberwithin{equation}{section}

\newtheorem{corollary}[theorem]{Corollary}

\newtheorem{definition}[theorem]{Definition}

\newtheorem{lemma}[theorem]{Lemma}

\newtheorem{proposition}[theorem]{Proposition}

\newenvironment{proof}[1][Proof]{\textbf{#1.}}{\ \rule{0.5em}{0.5em}}%

\begin{document}
\parindent 9mm
\title{Riemann-Liouville Fractional Cosine Functions
\thanks{This work was supported by the Natural Science Foundation of China (grant nos. 11301412 and 11131006), Research Fund for the Doctoral Program of Higher Education of China (grant no. 20130201120053), Shaanxi Province Natural Science Foundation of China (grant no. 2014JQ1017), Project funded by China Postdoctoral Science Foundation (grant no. 2014M550482), the Fundamental Research Funds for the Central Universities (grant no. 2012jdhz52).} \thanks{{\bf Mathematics Subject
Classification 2010} 33E12}}
\author{Zhan-Dong Mei \thanks{Corresponding
author, School of Mathematics and Statistics, Xi'an Jiaotong
University, Xi'an 710049, China; Email: zhdmei@mail.xjtu.edu.cn } \
\ \
 Ji-Gen Peng \thanks{School of Mathematics and Statistics, Xi'an Jiaotong
University, Xi'an 710049, China; Email: jgpeng@mail.xjtu.edu.cn}}


\date{}
\maketitle \thispagestyle{empty}
\begin{abstract}
%
%

In this paper, a new notion, named Riemann-Liouville fractional cosine function is presented. It is proved that a Riemann-Liouville $\alpha$-order fractional cosine function is equivalent to Riemann-Liouville $\alpha$-order fractional
resolvents introduced in [Z.D. Mei, J.G. Peng, Y. Zhang, Math. Nachr. 288, No. 7, 784-797 (2015)].

\vspace{0.5cm} 

%
%
\noindent {\bf Key words:} Riemann-Liouville fractional cosine function;
fractional resolvents; fractional differential equations.

\end{abstract}


\section{Introduction}

Assume that $X$ is a
Banach space, $A: D(A)\subset X\rightarrow X$ and $B: D(B)\subset X\rightarrow X$ are closed linear
operators. It is well-known that $C_0$-semigroup is an important
tool to study the following abstract Cauchy problem of first order
\begin{align}\label{firstorder}
    \left\{
      \begin{array}{ll}
        \frac{du(t)}{dt}=Au(t), & \hbox{} t>0\\
        u(0)=x, & \hbox{ }
      \end{array}
    \right.
\end{align}
and cosine function essentially characterizes the abstract Cauchy problem of second order described by
\begin{align}\label{secondorder}
    \left\{
      \begin{array}{ll}
        \frac{d^2u(t)}{dt^2}=Bu(t), & \hbox{} t>0\\
        u(0)=x, u'(0)=0. & \hbox{ }
      \end{array}
    \right.
\end{align}
Here a $C_0$-semigroup is a family $\{T(t)\}_{t\geq 0}$ of strongly continuous and bounded linear operators defined on $X$ satisfying $T(0)=I$ and $T(t+s)=T(t)T(s),\ t,s\geq 0$; a cosine function is a family $\{S(t)\}_{t\geq 0}$ of strongly continuous and bounded linear operators defined on $X$ satisfying $S(0)=I$ and $2S(t)S(s)=S(t)+S(s),\ t\geq s\geq 0.$

Concretely, system (\ref{firstorder}) is well-posed if and only if $A$ generates a $C_0$-semigroup $\{T(t)\}_{t\geq 0}$, namely, $Ax=\lim_{t\rightarrow 0^+}t^{-1}(T(t)x-x)$ with domain $D(A)=\{x\in D(A):\lim_{t\rightarrow 0^+}t^{-1}(T(t)x-x)$ exists $\}$;
system (\ref{secondorder}) is well-posed if and only if $B$ generates a cosine function $\{S(t)\}_{t\geq 0}$, namely, $Bx=2\lim_{t\rightarrow 0^+}t^{-2}(S(t)x-x)$ with domain $D(B)=\{x\in D(B):\lim_{t\rightarrow 0^+}t^{-2}(S(t)x-x)$ exists $\}$.
Therefore, pure algebraic methods can be used to study
abstract Cauchy problems of first and second orders.
For details, we refer to \cite{Engel2000,Goldstein1985}.

However, equations of integer order such as (\ref{firstorder}) and (\ref{secondorder}) can't exactly describe the behavior of many physical systems; fractional differential equations maybe more suitable for describing anomalous diffusion
on fractals (physical objects of fractional dimension, like some amorphous semiconductors
or strongly porous materials; see \cite{Anh2001,Metzler2001} and the references therein), fractional random walk \cite{Germano2009,Scalas2006},
etc. Fractional derivatives appear in the theory
of fractional differential equations; they describe the property of
memory and heredity of materials, and it is the major advantage of
fractional derivatives compared with integer order derivatives. Let $\alpha>0$ and $m=[\alpha]$, The smallest integer larger than or equal to $\alpha$. There are mainly two types of $\alpha$-order fractional differential equations, which are most used in the real problems.

1) Caupto fractional abstract Cauchy problem
\begin{align}\label{caputo}
    \left\{
      \begin{array}{ll}
        ^CD_t^\alpha u(t)=Au(t), & \hbox{ } t>0,\\
        u(0)=x, u^{(k)}(0)=0,k=1,2, ..., m-1.
      \end{array}
    \right.
\end{align}
where $^CD_t^\alpha$ is the Caupto fractional differential operator
defined as follows:
\begin{align*}
    ^CD_t^\alpha
    u(t)=\frac{1}{\Gamma(m-\alpha)}\int_0^t(t-\sigma)^{-\alpha}u^{(m)}(\sigma)d\sigma;
\end{align*}

 2) Riemann-Liouville fractional abstract Cauchy
problem
\begin{align}\label{R-L}
    \left\{
      \begin{array}{ll}
        D_t^\alpha u(t)=Au(t), & \hbox{ } \\
        (g_{2-\alpha}*u)(0)=\lim_{s\rightarrow 0^+}\int_0^s\frac{(s-\sigma)^{m-1-\alpha}}{\Gamma(2-\alpha)}u(\sigma)d\sigma= x, & \hbox{ } \\
 (g_{2-\alpha}*u)^{(k)}(0)=\lim_{s\rightarrow 0^+}\int_0^s\frac{d^k}{dt^k}\frac{(s-\sigma)^{m-1-\alpha}}{\Gamma(m-\alpha)}u(\sigma)d\sigma
=0, k=1,2, ..., m-1.& \hbox{ }
      \end{array}
    \right.
\end{align}
where $D_t^\alpha$ is the Riemann-Liouville fractional differential
operator defined by
\begin{align*}
    D_t^\alpha
    u(t)=\frac{1}{\Gamma(m-\alpha)}\frac{d}{dt}\int_0^t(t-\sigma)^{m-1-\alpha}u(\sigma)d\sigma.
\end{align*}

Obviously, (\ref{firstorder}) is just the limit state of equations (\ref{caputo}) and (\ref{R-L}) as $\alpha\rightarrow 1$, and (\ref{secondorder}) is just the limit state of equations (\ref{caputo}) and (\ref{R-L}) as $\alpha\rightarrow 2$. Initial conditions for the Caputo fractional derivatives are expressed in terms of initials of integer order derivatives \cite{Eidelman2004,Mei2014,Peng2012}. For some real materials, initial conditions should be expressed in terms of
Riemann-Liouville fractional derivatives, and it is possible to
obtain initial values for such initial conditions by appropriate
measurements \cite{Heymans2006, Hilfer2000}.

In order to study Caputo fractional abstract Cauchy problem  (\ref{caputo}), Bajlekova \cite{Bazhlekova2001} introduced the important notion of solution operator for equations (\ref{caputo}) as follows.

\begin{definition}\label{solutiona}
A family $\{T(t)\}_{t\geq 0}$ of bounded linear operators of $X$ is
called a solution operator for (\ref{caputo}) if the following three
conditions are satisfied:

(a) $T(t)$ is strongly continuous for $t\geq 0$ and $T(0)=I$,

(b) $T(t)D(A)\subset D(A)$ and $AT(t)x=T(t)Ax$ for all $x\in D(A)$
and $t\geq 0$,

(c) for any $x\in D(A)$, there holds
\begin{align*}
    T(t)x=x+J_t^\alpha T(t)Ax,
    t\geq 0.
\end{align*}
Here the notation $J_t^\alpha f(t)$ is defined by
\begin{align*}
    J_t^\alpha
f(t)=\frac{1}{\Gamma(\alpha)}\int_0^t(t-\sigma)^{\alpha-1}f(t)dt.
\end{align*}
\end{definition}

Chen and Li developed in \cite{Chen2010} a notion of $\alpha$-resolvent operator function, which was proved to be a new characteristic of solution operator. Hence, Caputo fractional abstract Cauchy problem can be studied by pure algebraic methods.
The definition of $\alpha$-resolvent operator function is as follows.
\begin{definition}\label{lichen}
Let $\{S(t)\}_{t\geq 0}$ be a family of bounded linear operators on
$X$. Then $\{S(t)\}_{t\geq 0}$  is called to be an $\alpha$-resolvent operator function, if the following assumptions are satisfied:

1) $S(t)$ is strongly continuous and $S(0)=I.$

2) $S(s)S(t)=S(t)S(s)$ for all $t,s\geq 0$.

3) $S(s)J_t^\alpha S(t)-J_s^\alpha S(s)S(t)=J_t^\alpha
S(t)-J_s^\alpha S(s)$ for all $t,s\geq 0$.
\end{definition}

In \cite{Li2012}, Li and Peng proposed the following notion of fractional resolvent to study Riemann-Liouville $\alpha$-order fractional abstract Cauchy problem (\ref{R-L}) with $\alpha\in (0,1)$.
\begin{definition}\cite{Li2012}\label{de}
 Let $0<\alpha< 1$. A family $\{T(t)\}_{t>0}$ of bounded linear operators on Banach space $X$ is called an
 $\alpha$-order fractional resolvent
 if it satisfies the following assumptions:

(P1) for any $x\in X$, $\ T(\cdot)x\in C((0,\infty),X)$, and
\begin{equation}\label{clear}
\lim_{t\rightarrow 0+}\Gamma(\alpha)t^{1-\alpha}T(t)x=x \ \
\mbox{for all } \ x\in X;
\end{equation}

(P2) $T(s)T(t)=T(t)T(s) \ \ \mbox{for all } t,s>0;$

(P3) for all $t,s>0$, there holds
\begin{equation}\label{fi}
\ T(t) J_{s}^{\alpha}T(s)-
J_{t}^{\alpha}T(t)T(s)=\frac{t^{\alpha-1}}{\Gamma(\alpha)}J_{s}^{\alpha}T(s)
-\frac{s^{\alpha-1}}{\Gamma(\alpha)}J_{t}^{\alpha}T(t).
\end{equation}
\end{definition}

Recently, we studied in \cite{Mei2015} Riemann-Liouville $\alpha$-order fractional Cauchy problem (\ref{R-L}) with order $\alpha\in (1,2)$ through the study of Riemann-Liouville
$\alpha$-order fractional resolvent defined as follows:
\begin{definition}\label{de}
 A family $\{T(t)\}_{t>
0}$ of bounded linear operators is called Riemann-Liouviille
$\alpha$-order fractional resolvent if it satisfies the following
assumptions:

(a) For any $x\in X$, $\ T_{\alpha}(\cdot)x\in C((0,\infty),X)$, and
\begin{equation}\label{clear}
\lim_{t\rightarrow 0^+}\Gamma(\alpha-1)t^{2-\alpha}T(t)x=x \ \
\mbox{for all } \ x\in X;
\end{equation}

(b) $ \ T(s)T_{\alpha}(t)=T(t)T_{\alpha}(s) \ \
\mbox{for all } t,s>0;$

(c) for all $t,s>0$, there holds
\begin{equation}\label{fi}
\ T(s) J_{t}^{\alpha}T(t)-
J_{s}^{\alpha}T(s)T(t)=\frac{s^{\alpha-2}}{\Gamma(\alpha-1)}J_{t}^{\alpha}T(t)
-\frac{t^{\alpha-2}}{\Gamma(\alpha-1)}J_{s}^{\alpha}T(s).
\end{equation}
\end{definition}

The linear operator $A$ defined by
\begin{align*}
D(A)=\{x\in X: \lim_{t\rightarrow
0^+}\frac{t^{1-\alpha}T(t)x-\frac{1}{\Gamma(\alpha)}x}{t^{2\alpha}}\
\mbox {exists}\}
\end{align*}
and
\begin{align*}
Ax=\lim_{t\rightarrow
0^+}\frac{t^{1-\alpha}T(t)x-\frac{1}{\Gamma(\alpha)}x}{t^{2\alpha}}
\ \ \mbox {for}\ x \in D(A)
\end{align*}
is the generator of the Riemann-Liouville $\alpha$-order fractional resolvent $\{T(t)\}_{t>0}$ in Definition \ref{de} with $D(A)$ the domain of $A$.

Moreover, we proved that $\{T(t)\}_{t>0}$ is a Riemann-Liouville $\alpha$-order fractional resolvent if and only if it is a solution operator defined as follows:

\begin{definition}\label{solutiona}
A family $\{T(t)\}_{t>0}$ of bounded linear operators of $X$ is
called a solution operator for (\ref{R-L}) if the following three
conditions are satisfied:

(a) $T(t)$ is strongly continuous for $t> 0$ and $\lim_{t\rightarrow 0^+}\Gamma(\alpha-1)t^{2-\alpha}T(t)x=x,\ x\in X$,

(b) $T(t)D(A)\subset D(A)$ and $AT(t)x=T(t)Ax$ for all $x\in D(A)$
and $t> 0$,

(c) for any $x\in D(A)$, there holds
\begin{align*}
    T(t)x=\frac{t^{1-\alpha}}{\Gamma(2-\alpha)}x+J_t^\alpha T(t)Ax,
    t> 0.
\end{align*}
\end{definition}

However, the above functional equations for fractional differential equations are not expressed in terms of the sum of time variables: $s+t$. This is very important in concrete
applications of the functional equation, just like $C_0$-semigroups,
cosine functions. Motivated by this, Peng and Li \cite{Peng2012} established the characteristic of $\alpha$-order fractional semigroup with $\alpha\in (0,1)$:

\begin{align*}
   \nonumber &\int_0^{t+s}\frac{T(\tau)}{(t+s-\tau)^\alpha}d\tau-
    \int_0^t\frac{T(\tau)}{(t+s-\tau)^\alpha}d\tau-
    \int_0^s\frac{T(\tau)}{(t+s-\tau)^\alpha}d\tau\\
    =&\alpha\int_0^t\int_0^s\frac{T(r_1)T(r_2)}{(t+s-r_1-r_2)^{1+\alpha}}dr_1dr_2,t,s\geq 0,
\end{align*}
where the integrals are in the sense of strong operator topology.
Concretely, they proved that $\alpha$-order
fractional semigroup is closely related to the solution operator of
Caputo fractional abstract Cauchy problem (\ref{caputo}).

Mei, Peng and Zhang \cite{Mei2013} developed the notion of Riemann-Liouville fractional semigroup as follows:
\begin{definition}
We call a family $\{T(t)\}_{t>0}$ of bounded linear operators to be
a Riemann-Liouville $\alpha$-order fractional semigroup on Banach
space $X$, if the following conditions are satisfied:

i) for any $x\in X$, $t\mapsto T(t)x$ is continuous over
$(0,\infty)$ and
\begin{equation}\label{clear}
\lim_{t\rightarrow 0+}\Gamma(\alpha)t^{1-\alpha}T(t)x=x;
\end{equation}

ii) for all $t,s> 0$, there holds
\begin{align}\label{cosin}
\Gamma(1-\alpha)T(t+s)=\alpha\int_0^t\int_0^s
\frac{T(r_1)T(r_2)}{(t+s-r_1-r_2)^{1+\alpha}}dr_1dr_2,
\end{align}
where the integrals are in the sense of strong operator topology.
\end{definition}
It is proved in \cite{Mei2013} that $A$ generates a Rimann-Liouville fractional semigroup if and only if it generates a fractional resolvent developed in \cite{Li2010}.

In order to study Caputo fractional Cauchy problem of order $\alpha\in (1,2)$, we recently studied in \cite{Mei2014} the notion of fractional cosine function as follows
\begin{definition}
A family $\{T(t)\}_{t\geq 0}$ of bounded and strongly continuous operators is called an $\alpha$-fractional cosine function if $T(0)=I$ and there hold
\begin{align}\label{cosin}
    \nonumber&\int_0^{t+s}\int_0^\sigma\frac{T(\tau)}{(t+s-\sigma)^{\alpha-1}}d\tau
    d\sigma-\int_0^t\int_0^\sigma\frac{T(\tau)}{(t+s-\sigma)^{\alpha-1}}d\tau
    d\sigma\\
    \nonumber &-\int_0^s\int_0^\sigma\frac{T(\tau)}{(t+s-\sigma)^{\alpha-1}}d\tau
    d\sigma\\
    \nonumber =&\int_0^t\int_0^s\frac{T(\sigma)T(\tau)}{(t-\sigma)^{\alpha-1}}d\tau
    d\sigma+\int_0^t\int_0^s\frac{T(\sigma)T(\tau)}{(s-\tau)^{\alpha-1}}d\tau
    d\sigma\\
    &-\int_0^t\int_0^s\frac{T(\sigma)T(\tau)}{(t+s-\sigma-\tau)^{\alpha-1}}d\tau
d\sigma, \ t,s\geq 0,
\end{align}
where the integrals are in the sense of strong operator topology.
\end{definition}
We proved that $A$ generates a fractional cosine function $\{T(t)\}_{t\geq 0}$ if and only if it generates an $\alpha$-resolvent operator function, that is, the following equalities holds:
\begin{equation*}
T(s)J_t^\alpha T(t)-J_s^\alpha T(s)T(t)=J_t^\alpha
T(t)-J_s^\alpha T(s), \ t,s\geq 0.
\end{equation*}

As is stated above, functional equations involving $t$, $s$ and $t+s$ have been discussed for Caputo fractional differential equations (\ref{caputo}) with $\alpha\in (0,1)$ and $\alpha\in (1,2)$,  Riemann-Liouville fractional equation (\ref{R-L}) with $\alpha\in (0,1)$. To close the gap, we will discuss the residual case, that is, functional equations involving $t$, $s$ and $t+s$ for Riemann-Liouville fractional equation (\ref{R-L}) with $\alpha\in (1,2)$. To this end, we first consider the special case that $T(\cdot)$ is exponentially
bounded (hence it is Laplace transformable).
Take laplace transform on both sides of (\ref{fi}) with respect to $s$ and $t$ to get
\begin{align}\label{lap}
    (\lambda^{-\alpha}-\mu^{-\alpha})\hat{T}(\mu)\hat{T}(\lambda)
    =\lambda^{1-\alpha}\mu^{1-\alpha}(\lambda^{-1}\hat{T}(\lambda)
    -\mu^{-1}\hat{T}(\mu)).
\end{align}
It follows from \cite[(3.8)]{Mei2014} that the Laplace transform of the right side of (\ref{cosin}) satisfies
\begin{align}\label{cosin1}
    \nonumber&\int_0^\infty e^{-\mu t}\int_0^\infty e^{-\lambda s}\bigg(\int_0^t\int_0^s\frac{T(\sigma)T(\tau)}{(t-\sigma)^{\alpha-1}}d\tau
    d\sigma+\int_0^t\int_0^s\frac{T(\sigma)T(\tau)}{(s-\tau)^{\alpha-1}}d\tau
    d\sigma\\
    \nonumber &-\int_0^t\int_0^s\frac{T(\sigma)T(\tau)}{(t+s-\sigma-\tau)^{\alpha-1}}d\tau
d\sigma\bigg) dsdt\\
=&\frac{\Gamma(2-\alpha)(\lambda^\alpha-\mu^\alpha)}{\lambda\mu(\lambda-\mu)}
\hat{T}(\mu)\hat{T}(\lambda).
\end{align}
The combination of (\ref{lap}) and (\ref{cosin1}) implies that
\begin{align*}
    &\int_0^\infty e^{-\mu t}\int_0^\infty e^{-\lambda s}\bigg(\int_0^t\int_0^s\frac{T(\sigma)T(\tau)}{(t-\sigma)^{\alpha-1}}d\tau
    d\sigma+\int_0^t\int_0^s\frac{T(\sigma)T(\tau)}{(s-\tau)^{\alpha-1}}d\tau
    d\sigma\\
    &-\int_0^t\int_0^s\frac{T(\sigma)T(\tau)}{(t+s-\sigma-\tau)^{\alpha-1}}d\tau
d\sigma\bigg) dsdt\\
=&\frac{\Gamma(2-\alpha)(\lambda^{-1}\hat{T}(\lambda)
-\mu^{-1}\hat{T}(\mu))}{\mu-\lambda}.
\end{align*}
Let $m(t)=\int_0^t T(\sigma)d\sigma$, by similar proof of \cite[
(4.2)]{Keyantuo2009}, there holds
\begin{align*}
    \int_0^\infty e^{-\mu t}\int_0^\infty e^{-\lambda s}
m(t+s)dsdt=\frac{\hat{m}(\mu)-\hat{m}(\lambda)}{\lambda-\mu}
=\frac{\lambda^{-1}\hat{T}(\lambda)
-\mu^{-1}\hat{T}(\mu)}{\mu-\lambda}.
\end{align*}
By virtue of Laplace transform, it follows that
\begin{align}\label{cosin2}
    \nonumber&\Gamma(2-\alpha)\int_0^{t+s}T(\sigma)d\sigma\\
    \nonumber =&\int_0^t\int_0^s\frac{T(\sigma)T(\tau)}{(t-\sigma)^{\alpha-1}}d\tau
    d\sigma+\int_0^t\int_0^s\frac{T(\sigma)T(\tau)}{(s-\tau)^{\alpha-1}}d\tau
    d\sigma\\
    &-\int_0^t\int_0^s\frac{T(\sigma)T(\tau)}{(t+s-\sigma-\tau)^{\alpha-1}}d\tau
d\sigma,
\end{align}

In the following two sections, we will show that (\ref{cosin2}) also holds without the assumption that $\{T(t)\}_{t>0}$ is exponentially bounded and it essentiality describes a Riemann-Liouville fractional resolvent.

\section{Riemann-Liouviller Fractional Cosine Function}

Equality (\ref{fi}) is an important functional equation for the solution of equation (\ref{R-L}) with $\alpha\in(1,2)$. However, as is stated in the introduction, (\ref{fi}) does not write the functional equation in terms of the sum of time variables: $s+t$. This is very important in concrete applications of the algebraic functional equation. Therefore, it is very valuable to study  functional equation (\ref{cosin2}), which appears in the following definitions.

\begin{definition}\label{DEF}
We call a family $\{T(t)\}_{t\geq 0}$ of bounded linear operators to
be a Riemann-Liouville $\alpha$-order fractional cosine function on Banach space $X$, if the following conditions are satisfied:

i) $T(t)$ is strongly continuous, that is, for any $x\in X$, the
mapping $t\mapsto T(t)x$ is continuous over $(0,\infty)$;

ii) there holds that
\begin{equation}\label{clear}
\lim_{t\rightarrow
0+}t^{2-\alpha} T(t)x
=\frac{x}{\Gamma(\alpha-1)} \ \ \mbox{for all } \ x\in X;
\end{equation}

iii) for all $t,s> 0$, there holds
\begin{align}\label{sin}
    \nonumber&
    \Gamma(2-\alpha)\int_0^{t+s}T(\sigma)d\sigma\\
    \nonumber =&\int_0^t\int_0^s\frac{T(\sigma)T(\tau)}{(t-\sigma)^{\alpha-1}}d\tau
    d\sigma+\int_0^t\int_0^s\frac{T(\sigma)T(\tau)}{(s-\tau)^{\alpha-1}}d\tau
    d\sigma\\
    &-\int_0^t\int_0^s\frac{T(\sigma)T(\tau)}{(t+s-\sigma-\tau)^{\alpha-1}}d\tau
d\sigma,
\end{align}
where the integrals are in the sense of strong operator topology.
\end{definition}

\begin{lemma}\label{commutative}
Let $\{T(t)\}_{t> 0}$ be a Riemann-Liouville $\alpha$-order fractional cosine on Banach space $X$. Then $\{T(t)\}_{t> 0}$ is commutative, i.e. $T(t)T(s)=T(s)T(t)$ for all $t,s > 0$.
\end{lemma}

\begin{proof}\ \
Observe that the left side of (\ref{sin}) is symmetric with respect to $t$ and $s$. Hence we can obtain the following equality.
\begin{align*}
&\int_0^t\int_0^s\frac{T(\sigma)T(\tau)}{(t-\sigma)^{\alpha-1}}d\tau
    d\sigma+\int_0^t\int_0^s\frac{T(\sigma)T(\tau)}{(s-\tau)^{\alpha-1}}d\tau
    d\sigma
    -\int_0^t\int_0^s\frac{T(\sigma)T(\tau)}{(t+s-\sigma-\tau)^{\alpha-1}}d\tau
d\sigma\\
=&
\int_0^s\int_0^t\frac{T(\sigma)T(\tau)}{(s-\sigma)^{\alpha-1}}d\tau
    d\sigma+\int_0^s\int_0^t\frac{T(\sigma)T(\tau)}{(t-\tau)^{\alpha-1}}d\tau
    d\sigma
    -\int_0^s\int_0^t\frac{T(\sigma)T(\tau)}{(t+s-\sigma-\tau)^{\alpha-1}}d\tau
d\sigma,\ t,s> 0.
\end{align*}
The commutative is proved by the same procedure of
\cite[Proposition 3.4]{Mei2013}.
\end{proof}

\begin{definition}\label{sss}
Let $\{T(t)\}_{t>0}$ be a Riemann-Liouville $\alpha$-order fractional cosine function on Banach space $X$.
Denote by $D(A)$ the set of all $x\in X$ such that the limit
\begin{align*}
    \lim_{t\rightarrow 0^+}\Gamma(\alpha+1)t^{-\alpha}J_t^{2-\alpha}\bigg(T(t)x-
    \frac{t^{\alpha-2}}{\Gamma(\alpha-1)}x\bigg)
\end{align*}
exists. Then, the operator $A:D(A)\rightarrow X$ defined by
\begin{align*}
    Ax=\lim_{t\rightarrow 0^+}\Gamma(\alpha+1)t^{-\alpha}J_t^{2-\alpha}\bigg(T(t)x-
    \frac{t^{\alpha-2}}{\Gamma(\alpha-1)}x\bigg)
\end{align*}
is called the generator of $\{T(t)\}_{t> 0}.$
\end{definition}

\begin{proposition}\label{pr2}
Assume $\{T(t)\}_{t>0}$ to be a Riemann-Liouville $\alpha$-order fractional cosine function on Banach space $X$. Suppose that $A$ is the generator of $\{T(t)\}_{t>0}$. Then,

(a) For any $x\in X$ and $t>0$, there holds $J_t^\alpha T(t)x\in
D(A)$
 and
 \begin{align}\label{aabbcc}
    T(t)x=\frac{t^{\alpha-2}}{\Gamma(\alpha-1)}x+AJ_t^\alpha T(t)x;
 \end{align}

(b) $T(t)D(A)\subset D(A)$ and $T(t)Ax=AT(t)x$, for all $x\in
  D(A)$.

 (c) For all $x\in D(A)$, we have
\begin{align*}
    T(t)x=\frac{t^{\alpha-2}}{\Gamma(\alpha-1)}x+J_t^\alpha T(t)Ax;
 \end{align*}

(d) $A$ is equivalently defined by
\begin{align}\label{ssf}
    Ax=\Gamma(2\alpha-1)\lim_{t\rightarrow
0^+}\frac{T(t)x-\frac{t^{\alpha-2}}
{\Gamma(\alpha-1)}x}{t^{2\alpha-2}}
\end{align}
and $D(A)$ is just consists of those $x\in X$ such that the above
limit exists.

(e) $A$ is closed and densely defined.

(f) $A$ admits at most one Riemann-Liouville $\alpha$-order fractional cosine function.

\end{proposition}
\begin{proof}\ \
(a) Let $x\in X$ and $b>0$ be fixed. Denote by $g_{b}(\cdot)$ the
truncation of $T(\cdot)$ at $b$, that is,
\begin{align*}
 g_{b}(\sigma)=\left\{
                 \begin{array}{ll}
                   T(\sigma), & \hbox{ } 0< \sigma\leq b\\
                   0, & \hbox{} \sigma>b.
                 \end{array}
               \right.
\end{align*}
Define the function
$H_b(r,s)$ for $r,s> 0$ by
\begin{align}\label{Htrs}
    H_b(r,s)=\bigg(g_b(r)-\frac{r^{\alpha-2}}{\Gamma(\alpha-1)}I\bigg)J_s^\alpha
g_b(s)x.
\end{align}
Obviously, for $0< r\leq t$,
\begin{align}\label{HET}
    H_t(r,t)=\bigg(T(r)-\frac{r^{\alpha-2}}{\Gamma(\alpha-1)}I\bigg)J_t^\alpha T(t)x.
\end{align}
Take Laplace transform with respect to $r$ and $s$ successively for
both sides of (\ref{Htrs}) to derive
\begin{align}\label{htl}
    \hat{H}_b(\mu,\lambda)=\lambda^{-\alpha}\hat{g}_b(\mu)\hat{g}_b(\lambda)x
    -\lambda^{-\alpha}\mu^{1-\alpha}\hat{g}_b(\lambda)x.
\end{align}

Denote by $L(t,s)$ and $R(t,s)$ the left and right sides of equality
(\ref{sin}), respectively. Moreover, denote by $R_b(t,s)$, and $L_b(t,s)$ the quantities resulted by replacing
$T(t)$ with $g_b(t)$ in $R(t,s)$, $L(t,s)$,
respectively.

It follows from (3.7) of \cite{Mei2014} that the Laplace transform of $R_b(t,s)$ with respect to $t$ and $s$ is given by
\begin{align}\label{R}
\hat{R}_b(\mu,\lambda)
=\frac{\Gamma(2-\alpha)(\lambda^\alpha-\mu^\alpha)}{\lambda\mu(\lambda-\mu)}\hat{g}_b(\mu)\hat{g}_b(\lambda).
\end{align}

For all $t> 0$, the Laplace transform of $\hat{L}_b(t,s)$ with respect to $s$ and $t$ can be obtained as follows:
\begin{align}\label{L}
    \hat{L}_b(\mu,\lambda)
    =\Gamma(2-\alpha)\frac{\lambda^{-1}\hat{g_b}(\lambda)
-\mu^{-1}\hat{g_b}(\mu)}{\mu-\lambda}.
\end{align}

Combine (\ref{htl}), (\ref{R}) and (\ref{L}) to derive
\begin{align*}
  \hat{H}_b(\mu,\lambda)=&\mu^{-\alpha}\hat{g}_b(\mu)\hat{g}_b(\lambda)x
    -\mu^{-\alpha}\lambda^{1-\alpha}\hat{g}_b(\mu)x\\
    &+\frac{\lambda^{1-\alpha}\mu^{1-\alpha}(\lambda-\mu)}{\Gamma(2-\alpha)}(\hat{L}_b(\mu,\lambda)-\hat{R}_b(\mu,\lambda))x.
\end{align*}

Take inverse Laplace transform to obtain
\begin{align*}
  H_b(r,s)=&\bigg(g_b(s)-\frac{s^{\alpha-2}}{\Gamma(\alpha-1)}I\bigg)J_r^\alpha g_b(r)x\\
    &+\frac{[(D_s^{2-\alpha})J_r^{\alpha-1}-(D_r^{2-\alpha})J_s^{\alpha-1}]\cdot[L_b(r,s)-R_b(r,s)]x}{\Gamma(2-\alpha)}.
\end{align*}
Here the Laplace transform formulas
\begin{align*}
    \widehat{D_b^\beta f}(\lambda)=\lambda^\beta
    \hat{f}(\lambda)-\lim_{t\rightarrow 0^+}J_t^{\alpha-1}f(t), 0<\beta<1, f\in C([0,\infty),X)
\end{align*}
is used.

 By the definition of $g_b$, it follows that
$L_b(r,s)=R_b(r,s)$ for all $0< s,r\leq b$, we have that
\begin{align*}
    H_b(r,s)=\bigg(T(s)-\frac{s^{\alpha-2}}{\Gamma(\alpha-1)}I\bigg)J_r^\alpha T(r)x, \forall \mbox{ }0< r,s \leq b.
\end{align*}
This implies that
\begin{align}\label{HTT}
    H_t(r,t)=\bigg(T(t)-\frac{t^{\alpha-2}}{\Gamma(\alpha-1)}I\bigg)J_r^\alpha T(r)x,\forall \mbox{ } 0< r\leq t.
\end{align}
Combining (\ref{HET}) and (\ref{HTT}), we derive that
\begin{align*}
    &\lim_{r\rightarrow 0^+}\Gamma(\alpha+1)r^{-\alpha}J_r^{2-\alpha}\bigg(T(r)-\frac{r^{\alpha-2}}{\Gamma(\alpha-1)}\bigg)J_t^\alpha
    T(t)x\\
    =&\lim_{r\rightarrow
    0^+}\Gamma(\alpha+1)r^{-\alpha}\bigg(T(t)-\frac{t^{\alpha-2}}{\Gamma(\alpha-1)}I\bigg)J_r^2
    T(r)x\\
    =&\Gamma(\alpha+1)\bigg(T(t)-\frac{t^{\alpha-2}}{\Gamma(\alpha-1)}I\bigg)\lim_{r\rightarrow 0^+}r^{-\alpha}\int_0^r(r-\sigma)
    T(\sigma)xd\sigma\\
    =&\Gamma(\alpha+1)\bigg(T(t)-\frac{t^{\alpha-2}}{\Gamma(\alpha-1)}I\bigg)
\\
\cdot&\lim_{r\rightarrow
    0^+}\int_0^1(1-\sigma)\sigma^{\alpha-2}
    (r\sigma)^{2-\alpha}T(r\sigma)xd\sigma.
\end{align*}
By the dominated convergence theorem and (b) of Definition \ref{DEF}, it follows that
\begin{align*}
&\lim_{r\rightarrow 0^+}\Gamma(\alpha+1)r^{-\alpha}J_r^{2-\alpha}\bigg(T(r)-\frac{r^{\alpha-2}}{\Gamma(\alpha-1)}\bigg)J_t^\alpha
    T(t)x\\
    =&\Gamma(\alpha+1)\bigg(T(t)-\frac{t^{\alpha-2}}{\Gamma(\alpha-1)}I\bigg)
\int_0^1(1-\sigma)\sigma^{\alpha-2}
    \lim_{r\rightarrow
    0^+}(r\sigma)^{2-\alpha}T(r\sigma)xd\sigma\\
=&\frac{\Gamma(\alpha+1)}{\Gamma(\alpha-1)}\bigg(T(t)-\frac{t^{\alpha-2}}{\Gamma(\alpha-1)}I\bigg)
\int_0^1(1-\sigma)\sigma^{\alpha-2}d\sigma x\\
=&\frac{\Gamma(\alpha+1)}{\Gamma(\alpha-1)}\bigg(T(t)-\frac{t^{\alpha-2}}{\Gamma(\alpha-1)}I\bigg)
\frac{\Gamma(\alpha-1)\Gamma(2)}{\Gamma(\alpha+1)} x\\
=&T(t)x-\frac{t^{\alpha-2}}{\Gamma(\alpha-1)}x.
\end{align*}

This implies that $J_t^\alpha
    T(t)x\in D(A)$ and $$AJ_t^\alpha
    T(t)x=T(t)x-\frac{t^{\alpha-2}}{\Gamma(\alpha-1)}x.$$

(b) and (c) are directly obtained by Lemma \ref{commutative} and
(a).

(d) Denote by $D$ the set of those $x\in X$ such that the limit
$$\lim_{t\rightarrow
0^+}\frac{T(t)x-
\frac{t^{\alpha-2}}{\Gamma(\alpha-1)}x}{t^{2\alpha-2}}$$
exists. Let $x\in D(A).$ Then, by (b), we have that
\begin{align*}
    &\Gamma(2\alpha-1)\lim_{t\rightarrow
0^+}\frac{T(t)x-
\frac{t^{\alpha-2}}{\Gamma(\alpha-1)}x}{t^{2\alpha-2}}\\
    =&\Gamma(2\alpha-1)\lim_{t\rightarrow
    0^+}\frac{J_t^\alpha T(t)Ax}{t^{2\alpha-2}}\\
    =&\frac{\Gamma(2\alpha-1)}{\Gamma(\alpha)}\lim_{t\rightarrow
    0^+}\frac{1}{t^{2\alpha-2}}\int_0^t(t-\sigma)^{\alpha-1}T(\sigma)Axd\sigma\\
    =&\frac{\Gamma(2\alpha-1)}{\Gamma(\alpha)}
    \lim_{t\rightarrow
    0^+}\int_0^1(1-\sigma)^{\alpha-1}\sigma^{\alpha-2}
    (t\sigma)^{2-\alpha}T(t\sigma)Axd\sigma.
\end{align*}
The combination of the dominated convergence theorem and (b) of Definition \ref{DEF} indicates that
\begin{align*}
    &\Gamma(2\alpha-1)\lim_{t\rightarrow
0^+}\frac{T(t)x-
\frac{t^{\alpha-2}}{\Gamma(\alpha-1)}x}{t^{2\alpha-2}}\\
    =&\frac{\Gamma(2\alpha-1)}{\Gamma(\alpha)}
    \int_0^1(1-\sigma)^{\alpha-1}\sigma^{\alpha-2}
    \lim_{t\rightarrow
    0^+}(t\sigma)^{2-\alpha}T(t\sigma)Axd\sigma\\
    =&\frac{\Gamma(2\alpha-1)}{\Gamma(\alpha-1)\Gamma(\alpha)}
    \int_0^1(1-\sigma)^{\alpha-1}\sigma^{\alpha-2}Axd\sigma\\
    =&\frac{\Gamma(2\alpha-1)}{\Gamma(\alpha-1)\Gamma(\alpha)}
\frac{\Gamma(\alpha-1)\Gamma(\alpha)}{\Gamma(2\alpha-1)}Ax\\
    =&Ax.
\end{align*}
This implies that $x\in D$ and then $D(A)\subset D$. Now we prove
the converse inclusion.
Let $x\in D$, that is, the limit
$$\lim_{t\rightarrow
0^+}\frac{T(t)x-
\frac{t^{\alpha-2}}{\Gamma(\alpha-1))}x}{t^{2\alpha-2}}.$$
exists. By the dominated convergence theorem, it follows that
\begin{align*}
    &\lim_{t\rightarrow 0^+}\Gamma(\alpha+1)t^{-\alpha}J_t^{2-\alpha}\bigg(T(t)x-
    \frac{t^{\alpha-2}}{\Gamma(\alpha-1)}x\bigg)\\
    =&\lim_{t\rightarrow
    0^+}\frac{\Gamma(\alpha+1)}{\Gamma(2-\alpha)}
    \int_0^1(1-\sigma)^{1-\alpha}\sigma^{2\alpha-2}
   \frac{T(t\sigma)x-\frac{(t\sigma)^{\alpha-2}}{\Gamma(\alpha-1)}x}
{(t\sigma)^{2\alpha-2}}d\sigma\\
   =&\frac{\Gamma(\alpha+1)}{\Gamma(2-\alpha)}
    \int_0^1(1-\sigma)^{1-\alpha}\sigma^{2\alpha-2}
   \lim_{t\rightarrow
    0^+}\frac{T(t\sigma)x-\frac{(t\sigma)^{\alpha-2}}{\Gamma(\alpha-1)}x}
{(t\sigma)^{2\alpha-2}}d\sigma\\
    =&\frac{\Gamma(\alpha+1)}{\Gamma(2-\alpha)}
    \frac{\Gamma(2-\alpha)\Gamma(2\alpha-1)}{\Gamma(\alpha+1)} \lim_{t\rightarrow
    0^+}\frac{T(t)x-\frac{(t)^{\alpha-2}}{\Gamma(\alpha-1)}x}
{t^{2\alpha-2}}.
\end{align*}
Hence, $x\in D(A)$ and
\begin{align}\label{absc}
    Ax=\Gamma(2\alpha-1)\lim_{t\rightarrow
0^+}\frac{T(t)x-
\frac{t^{\alpha-2}}{\Gamma(\alpha-1))}x}{t^{2\alpha-2}}.
\end{align}

(e) The properties that $A$ is closed and densely defined are followed directly from the combination of (d) and \cite{Li2012}.

(f) Assume that both $\{T(t)\}_{t>0}$ and $\{S(t)\}_{t> 0}$
are Riemann-Liouville $\alpha$-order fractional resolvent generated
by $A$. Then, by (c), for all $x\in D(A)$, we have
\begin{align*}
    \frac{t^{\alpha-2}}{\Gamma(\alpha-1)}*T(t)x
=&(S(t)-J_t^\alpha AS(t))*T(t)x\\
=&S(t)*T(t)x-(J_t^\alpha AS(t))*T(t)x\\
=&S(t)*(T(t)x-J_t^\alpha AT(t)x)\\
=& \frac{t^{\alpha-2}}{\Gamma(\alpha-1)}*S(t)x.
 \end{align*}
By Titchmarsh's Theorem, for any $t> 0$, $T(t)=S(t)$ on $D(A)$. The result is obtained by the density of $A$.
\end{proof}

\begin{corollary}\label{zheng}
Assume that $A$ generates a Rimann-Liouville $\alpha$-order fractional cosine function on Banach space $X$. Then
$\{T(t)\}_{t> 0}$ is a Riemann-Liouville $\alpha$-order fractional resolvent.
\end{corollary}
\begin{proof}\ \
In (a) of Theorem \ref{pr2}, replacing $x$ with $J_s^\alpha T(s)x$,
and using Lemma \ref{commutative}, we obtain that
\begin{align*}
  T(t)J_s^\alpha T(s)x=&\frac{t^{\alpha-2}}{\Gamma(\alpha-1)}J_s^\alpha T(s)x
+AJ_t^\alpha T(t)J_s^\alpha T(s)x\\
=&\frac{t^{\alpha-2}}{\Gamma(\alpha-1)}J_s^\alpha T(s)x
+AJ_s^\alpha T(s)J_t^\alpha T(t)x\\
=&\frac{t^{\alpha-2}}{\Gamma(\alpha-1)}J_s^\alpha T(s)x
+\bigg(T(s)-\frac{t^{\alpha-2}}{\Gamma(\alpha-1)}\bigg)J_t^\alpha T(t)x,
\end{align*}
which is just (\ref{fi}). The proof is therefore completed.
\end{proof}

\section{Equivalent to Riemann-Liouville fractional resolvent}

In this section, we will prove that equality (\ref{fi}) essentially describes a Rimann-Liouville $\alpha$-order fractional cosine function.

\begin{theorem}\label{relation}
Suppose that $\{T(t)\}_{t> 0}$ is a Riemann-Liouville $\alpha$-order fractional resolvent on Banach space $X$.
Then, the family is a Riemann-Liouville $\alpha$-order fractional cosine function.
\end{theorem}

\begin{proof}\ \
Denote by $L(t,s)$ and $R(t,s)$ the left and right sides of equality
(\ref{sin}), respectively. Obviously, what we need is to prove
that $L(t,s)=R(t,s)$ for all $t,s> 0.$ For brevity, we introduce the
following notations. Let
\begin{align*}
    H(t,s)=&T(t)J_s^\alpha T(s)-J_t^\alpha T(t)T(s),\\
K(t,s)=&\frac{t^{\alpha-2}}{\Gamma(\alpha-1)}J_{s}^{\alpha}T(s)
-\frac{s^{\alpha-2}}{\Gamma(\alpha-1)}J_{t}^{\alpha}T(t), t,s> 0.
\end{align*}
Moreover, for sufficiently large $b>0$ denote by $g_b(t)$ the
truncation of $T(t)$ at $b$, and by $R_b(t,s)$, $L_b(t,s)$,
$H_b(t,s)$ and $K_b(t,s)$ the quantities resulted by replacing
$T(t)$ with $g_b(t)$ in $R(t,s)$, $L(t,s)$, $H(t,s)$ and $K(t,s)$,
respectively.

We set
\begin{align*}
    P_b(t,s)=&\int_0^t\int_0^s\frac{H_b(\sigma,\tau)}{(t-\sigma)^{\alpha-1}}d\tau
    d\sigma+\int_0^t\int_0^s\frac{H_b(\sigma,\tau)}{(s-\tau)^{\alpha-1}}d\tau
    d\sigma\\
    \nonumber &-\int_0^t\int_0^s\frac{H_b(\sigma,\tau)}{(t+s-\sigma-\tau)^{\alpha-1}}d\tau
d\sigma
\end{align*}
and
\begin{align}\label{limi}
    Q_b(t,s)=&\int_0^t\int_0^s\frac{K_b(\sigma,\tau)}{(t-\sigma)^{\alpha-1}}d\tau
    d\sigma+\int_0^t\int_0^s\frac{K_b(\sigma,\tau)}{(s-\tau)^{\alpha-1}}d\tau
    d\sigma\\
    \nonumber &-\int_0^t\int_0^s\frac{K_b(\sigma,\tau)}{(t+s-\sigma-\tau)^{\alpha-1}}d\tau
d\sigma.
\end{align}
Observe that the equality (\ref{fi}) implies $H(t,s)=K(t,s)$ for any $t,s> 0$. Thus, for all  $t,s>0$,
\begin{align}\label{li}
    \lim_{b\rightarrow \infty}P_b(t,s)=\lim_{b\rightarrow
\infty}Q_b(t,s).
\end{align}

By \cite[(3.13)]{Mei2014}, it follows that
\begin{align}\label{P}
    P_b(t,s)=(J_s^\alpha-J_t^\alpha)R_b(t,s), \forall\ t,s> 0.
\end{align}

We now compute Laplace transform of the first term of $Q_b(t,s)$
with respect to $s$ and $t$ as follows,
\begin{align*}
    &\int_0^\infty e^{-\mu t}\int_0^\infty e^{-\lambda s}
    \int_0^t\int_0^s\frac{K_b(\sigma,\tau)}{(t-\sigma)^{\alpha-1}}d\tau
    d\sigma dsdt\\
 =&\int_0^\infty e^{-\mu t}\int_0^\infty e^{-\lambda s}
    \int_0^t\int_0^s\frac{\frac{\sigma^{\alpha-2}}{\Gamma(\alpha-1)}
J_{\tau}^{\alpha}g_b(\tau)
-\frac{\tau^{\alpha-2}}{\Gamma(\alpha-1)}J_{\sigma}^{\alpha}g_b(\sigma)}{(t-\sigma)^{\alpha-1}}d\tau
    d\sigma dsdt\\
 =&\int_0^\infty e^{-\mu t}\int_0^t\int_0^\infty e^{-\lambda s}
    \int_0^s\frac{\frac{\sigma^{\alpha-2}}{\Gamma(\alpha-1)}J_{\tau}^{\alpha}g_b(\tau)
-\frac{\tau^{\alpha-2}}{\Gamma(\alpha-1)}J_{\sigma}^{\alpha}g_b(\sigma)}{(t-\sigma)^{\alpha-1}}d\tau
   ds d\sigma dt\\
 =&\int_0^\infty e^{-\mu t}\int_0^\infty e^{-\lambda s}
    \int_0^t\int_0^s\frac{\frac{\sigma^{\alpha-2}}{\Gamma(\alpha-1)}J_{\tau}^{\alpha}g_b(\tau)
-\frac{\tau^{\alpha-2}}{\Gamma(\alpha-1)}J_{\sigma}^{\alpha}g_b(\sigma)}{(t-\sigma)^{\alpha-1}}d\tau
    d\sigma dsdt\\
 =&\int_0^\infty e^{-\mu t}\int_0^t\frac{\frac{\sigma^{\alpha-2}}{\Gamma(\alpha-1)}
}{(t-\sigma)^{\alpha-1}}\int_0^\infty e^{-\lambda s}
    \int_0^sJ_{\tau}^{\alpha}g_b(\tau)d\tau
   ds d\sigma dt\\
 &-\int_0^\infty e^{-\mu t}\int_0^t\frac{J_{\sigma}^{\alpha}g_b(\sigma)
}{(t-\sigma)^{\alpha-1}}\int_0^\infty e^{-\lambda s}
    \int_0^s\frac{\tau^{\alpha-2}}{\Gamma(\alpha-1)}d\tau
   ds d\sigma dt\\
 =&\Gamma(2-\alpha)\mu^{-1}\lambda^{-\alpha-1}\hat{g}_b(\lambda)
 -\Gamma(2-\alpha)\mu^{-2}\lambda^{-\alpha}\hat{g}_b(\mu),
\end{align*}

The Laplace transform of the second term of $Q_b(t,s)$
with respect to $s$ and $t$ is computed as follows
\begin{align*}
    &\int_0^\infty e^{-\mu t}\int_0^\infty e^{-\lambda s}
    \int_0^t\int_0^s\frac{K_b(\sigma,\tau)}{(s-\tau)^{\alpha-1}}d\tau
    d\sigma dsdt\\
 =&\int_0^\infty e^{-\mu t}\int_0^\infty e^{-\lambda s}
    \int_0^t\int_0^s\frac{\frac{\sigma^{\alpha-2}}
{\Gamma(\alpha-1)}J_{\tau}^{\alpha}g_b(\tau)
-\frac{\tau^{\alpha-2}}{\Gamma(\alpha-1)}J_{\sigma}^{\alpha}g_b(\sigma)}{(s-\tau)^{\alpha-1}}d\tau
    d\sigma dsdt\\
 =&\int_0^\infty e^{-\mu t}\int_0^t\frac{\sigma^{\alpha-2}}{\Gamma(\alpha-1))}
\int_0^\infty e^{-\lambda s}
    \int_0^s\frac{J_{\tau}^{\alpha}g_b(\tau)
}{(s-\tau)^{\alpha-1}}d\tau
    dsd\sigma dt\\
 &-\int_0^\infty e^{-\mu t}\int_0^tJ_{\sigma}^{\alpha}g_b(\sigma)\int_0^\infty e^{-\lambda s}    \int_0^s\frac{\frac{\tau^{\alpha-2}}{\Gamma(\alpha-1)}
}{(s-\tau)^{\alpha-1}}d\tau
    dsd\sigma dt\\
 =&\Gamma(2-\alpha)\mu^{-\alpha}\lambda^{-2}\hat{g}_b(\lambda)-
 \Gamma(2-\alpha)\lambda^{-1}\mu^{-\alpha-1}\hat{g}_b(\mu).
\end{align*}
We compute the Laplace transform of the third term of $Q_b(t,s)$ with respect to $s$ and $t$ as follows.

\begin{align*}
    &-\int_0^\infty e^{-\mu t}\int_0^\infty e^{-\lambda s}
    \int_0^t\int_0^s\frac{K_b(\sigma,\tau)}{(t+s-\sigma-\tau)^{\alpha-1}}d\tau
d\sigma dsdt\\
=&-\int_0^\infty e^{-\mu t}\int_0^\infty e^{-\lambda s}
    \int_0^t\int_0^s\frac{\frac{\sigma^{\alpha-2}}
{\Gamma(\alpha-1)}J_{\tau}^{\alpha}g_b(\tau)
-\frac{\tau^{\alpha-2}}{\Gamma(\alpha-1)}J_{\sigma}^{\alpha}g_b(\sigma)}{(t+s-\sigma-\tau)^{\alpha-1}}d\tau
d\sigma dsdt\\
=&-\int_0^\infty e^{-\mu t}\int_0^t\frac{\sigma^{\alpha-2}}{\Gamma(\alpha-1)}\int_0^\infty e^{-\lambda s}
    \int_0^s\frac{J_{\tau}^{\alpha}g_b(\tau)
}{(t+s-\sigma-\tau)^{\alpha-1}}d\tau
 dsd\sigma dt\\
 &+\int_0^\infty e^{-\mu t}\int_0^tJ_{\sigma}^{\alpha}g_b(\sigma)\int_0^\infty e^{-\lambda s}    \int_0^s\frac{\frac{\tau^{\alpha-2}}{\Gamma(\alpha-1)}}{(t+s-\sigma-\tau)^{\alpha-1}}d\tau
 dsd\sigma dt\\
 =&-\int_0^\infty e^{-\mu t}\int_0^t\frac{\sigma^{\alpha-2}}{\Gamma(\alpha-1)}\int_0^\infty e^{-\lambda s}\frac{1}{(t+s-\sigma)^{\alpha-1}}
 dsd\sigma dt\lambda^{-\alpha}g_b(\lambda)\\
 &+\lambda^{1-\alpha}\int_0^\infty e^{-\mu t}\int_0^tJ_{\sigma}^{\alpha}g_b(\sigma)\int_0^\infty e^{-\lambda s}    \frac{1}{(t+s-\sigma)^{\alpha-1}}
 dsd\sigma dt\\
 =&-\int_0^\infty e^{-\mu t}\int_0^t\frac{\sigma^{\alpha-2}}{\Gamma(\alpha-1)}
 e^{\lambda(t-\sigma)}\bigg(\int_0^\infty e^{-\lambda r}r^{1-\alpha}
 dr-\int_0^{t-\sigma} e^{-\lambda r}r^{1-\alpha}
 dr\bigg)d\sigma dt\lambda^{-\alpha}g_b(\lambda)\\
 &+\lambda^{1-\alpha}\int_0^\infty e^{-\mu t}\int_0^tJ_{\sigma}^{\alpha}g_b(\sigma)e^{\lambda(t-\sigma)}
 \bigg(\int_0^\infty e^{-\lambda r} r^{1-\alpha}
 dr-\int_0^{t-\sigma} e^{-\lambda r} r^{1-\alpha}
 dr\bigg)d\sigma dt\\
 =&-\Gamma(2-\alpha)\lambda^{\alpha-2}
 \int_0^\infty e^{-\mu t}\int_0^t\frac{\sigma^{\alpha-2}}{\Gamma(\alpha-1)}
 e^{\lambda(t-\sigma)}d\sigma dt\lambda^{-\alpha}g_b(\lambda)\\
 &+\int_0^\infty e^{-\mu t}\int_0^t\frac{\sigma^{\alpha-2}}{\Gamma(\alpha-1)}
 \int_0^{t-\sigma} e^{\lambda(t-\sigma-r)}r^{1-\alpha}
 drd\sigma dt\lambda^{-\alpha}g_b(\lambda)\\
 &+\Gamma(2-\alpha)\lambda^{\alpha-2}\lambda^{1-\alpha}
 \int_0^\infty e^{-\mu t}\int_0^tJ_{\sigma}^{\alpha}g_b(\sigma)e^{\lambda(t-\sigma)}d\sigma dt\\
 &-\lambda^{1-\alpha}\int_0^\infty e^{-\mu t}\int_0^tJ_{\sigma}^{\alpha}g_b(\sigma)
 \int_0^{t-\sigma} e^{\lambda(t-\sigma-r)} r^{1-\alpha}
 drd\sigma dt\\
 =&-\Gamma(2-\alpha)\lambda^{\alpha-2}\lambda^{-\alpha}
 \frac{\mu^{1-\alpha}}{\mu-\lambda}g_b(\lambda)
 +\Gamma(2-\alpha)\frac{\mu^{-1}}{\mu-\lambda}\lambda^{-\alpha}g_b(\lambda)\\
 &+\Gamma(2-\alpha)\lambda^{\alpha-2}\lambda^{1-\alpha}
 \frac{\mu^{-\alpha}}{\mu-\lambda}\hat{g}_b(\mu)
 -\Gamma(2-\alpha)\lambda^{1-\alpha}\mu^{\alpha-2}
 \frac{\mu^{-\alpha}}{\mu-\lambda}\hat{g}_b(\mu).
\end{align*}
Using (\ref{L}), we can obtain that
\begin{align}\label{dd}
\nonumber\hat{Q}_b(\mu,\lambda)=&\int_0^\infty e^{-\mu t}\int_0^\infty e^{-\lambda s}
\bigg(\int_0^t\int_0^s\frac{K_b(\sigma,\tau)}{(t-\sigma)^{\alpha-1}}d\tau
d\sigma+\int_0^t\int_0^s\frac{K_b(\sigma,\tau)}{(s-\tau)^{\alpha-1}}d\tau
d\sigma\\
    \nonumber &-\int_0^t\int_0^s\frac{K_b(\sigma,\tau)}{(t+s-\sigma-\tau)^{\alpha-1}}d\tau
d\sigma)\bigg)dsdt\\
=&(\lambda^{-\alpha}-\mu^{-\alpha})\hat{L_b}(\mu,\lambda).
\end{align}
Taking inverse Laplace transform on both sides of (\ref{dd}), we derive
\begin{align}\label{Q}
    Q_b(t,s)=(J_s^\alpha-J_t^\alpha)L_b(t,s),\ \forall \ t, \ s> 0.
\end{align}
Form (\ref{P}) and (\ref{Q}), we have that
\begin{align*}
    (J_s^\alpha-J_t^\alpha)L(t,s)=(J_s^\alpha-J_t^\alpha)R(t,s),\ \forall \ t, \ s>0.
\end{align*}
Therefore, $L(t,s)=R(t,s)$. This completes the proof.
\end{proof}

Combining Corollary \ref{zheng} and Theorem \ref{relation},
we can obtain the equivalent of Riemann-Liouville $\alpha$-order fractional resolvents and Riemann-Liouville $\alpha$-order fractional cosine functions.

\end{document}